\documentclass[a4paper, 10pt, DIV10, headinclude=false, footinclude=false]{scrarticle}

\usepackage[utf8]{inputenc}
\usepackage[T1]{fontenc}
\usepackage[english]{babel}

\usepackage[
colorlinks=true,
urlcolor=blue,
linkcolor=black
]{hyperref}
\usepackage[pdftex]{graphicx}
\usepackage{subcaption}
\usepackage{latexsym}
\usepackage{amsmath,amssymb,amsthm}
\usepackage{esint} 
\usepackage{aligned-overset}
\usepackage{multirow}
\usepackage{bigdelim}
\usepackage{empheq} 
\usepackage{tikz}
\usetikzlibrary{positioning,shapes.multipart}
\newtheorem{Satz}{Theorem}[section]
\newtheorem{Lemma}[Satz]{Lemma}	
\newtheorem{Proposition}[Satz]{Proposition} 
\theoremstyle{definition}

\theoremstyle{remark}	 
\newtheorem{Bsp}[Satz]{Example}	  
\newtheorem{Bemerkung}[Satz]{Remark}                  
\numberwithin{equation}{section} 

\newenvironment{abstr}[1]{ \vspace{.05in}\footnotesize
	\parindent .2in
	{\upshape\bfseries #1. }\ignorespaces}{\par\vspace{.1in}}
\newenvironment{Abstract}{\begin{abstr}{Abstract}}{\end{abstr}}
\newenvironment{keywords}{\begin{abstr}{Key words}}{\end{abstr}}
\newenvironment{AMS}{\begin{abstr}{AMS subject classifications}}{\end{abstr}}

\newcommand{\R}{\mathbb{R}} 
\def\avint{\fint}
\def\a0n{\mathop{a_{n,K}^{\epsilon}}}
\def\chin{\mathop{\chi_n^{\epsilon}}}
\def\A0nK{A_0(t_n,x_K)}
\def\PhiH{\Phi_{H \mid_{I_{\epsilon,K}}}}
\def\phidisc{\phi_{h,k}^{\epsilon}}
\def\Bhh{B_{H,h}}
\def\Ahh{A_{H,h}}

\def\Vmicro{V_h^p}
\usepackage{lipsum}

\allowdisplaybreaks[4]

\begin{document}  

	\title{Fully discrete Heterogeneous Multiscale Method for parabolic problems with multiple spatial and temporal scales%
	\thanks{Funded by the Deutsche Forschungsgemeinschaft (DFG, German Research Foundation) under VE 1397/2-1. Major parts of this work were accomplished while BV was affiliated with Karlsruher Institut für Technologie.}
	}
	\author{Daniel Eckhardt\footnotemark[2]\and Barbara Verf\"urth\footnotemark[3]}
	\date{}
	\maketitle
	
	\renewcommand{\thefootnote}{\fnsymbol{footnote}}
	\footnotetext[2]{Institut für Angewandte und Numerische Mathematik, Karlsruher Institut für Technologie, Englerstr. 2, D-76131 Karlsruhe}
	\footnotetext[3]{Institut für Numerische Simulation, Universit\"at Bonn, Friedrich-Hirzebruch-Allee 7, D-53115 Bonn}
	\renewcommand{\thefootnote}{\arabic{footnote}}

	\begin{Abstract} The aim of this work is the numerical homogenization of a parabolic problem with several time and spatial scales using the heterogeneous multiscale method. 
	We replace the actual cell problem with an alternate one, using Dirichlet boundary and initial values instead of periodic boundary and time conditions. 
	Further, we give a detailed a priori error analysis of the fully discretized, i.e., in space and time for both the macroscopic and the cell problem, method. Numerical experiments illustrate the theoretical convergence rates.
	\end{Abstract}

    \begin{keywords}
    	multiscale method; numerical homogenization; parabolic problem; time-space multiscale coefficient;
    	a priori error estimates
    \end{keywords}
    
    \begin{AMS}
    	65M60, 65M15, 65M12, 35K15, 80M40	
    \end{AMS}

	\section{Introduction \label{Motivation}}
	Problems with multiple spatial and temporal scales occur in a variety  of different phenomena and materials.  Prominent examples are saltwater intrusion, storage of radioactive waste products or various composite materials (\cite{held2005homogenization,kamga2007numerical,kalamkarov2009asymptotic}). These examples all have in common that both macroscopic and microscopic scales occur. Consequently, they are particularly challenging from a numerical point of view. However, from the application point of view, it is often sufficient to know a description of the macroscopic properties. Therefore, it is quite relevant to develop a method that includes all small-scale effects without having to calculate them simultaneously. This
	is the main component of (numerical) homogenization. 
	
	In this work we are interested in the following parabolic problem 
	\begin{align*}
		\dfrac{\partial u^{\epsilon}}{\partial t} - \nabla \cdot \Bigl(a\Bigl(t,x, \frac{t}{\epsilon^2}, \frac{x}{\epsilon}\Bigr) \nabla u^{\epsilon}\Bigr)  &= f, 
	\end{align*}
	with initial and boundary condtions. The precise setting is given further below. $a\bigl(t,x, \frac{t}{\epsilon^2}, \frac{x}{\epsilon}\bigr)$ is called the time-space  multiscale coefficient and represents physical properties of the considered material. 
	If we use standard finite element  and time stepping methods, we obtain sufficiently good solutions only for small time steps \emph{and} fine grids as the following example illustrates. 
	\begin{Bsp}
		\label{Bsp}
		Let $\Omega = (0,1)$ and $T = 1$. Furthermore we consider
		\begin{align*}
			a(t,x,s,y) = 3 + \cos(2\pi y) + \cos^2(2\pi s).
		\end{align*}
		Let the initial condition be $ u^{\epsilon}(0,x) = 0$ for all $x \in \Omega$. 
		Figure \ref{Model_L2} shows the error in the $L^2$-norm with respect to the numerical and a reference solution.  The reference solutions were calculated using finite elements with grid width $h= 10^{-6}$ and the implicit Euler method time step size $\tau = 1/100$. The theory yields an expected quadratic order of convergence. However, this occurs here only for small grid sizes. 
		More precisely, the error converges only when $h < \epsilon$, see Figure \ref{Model_L2}. 
		Similar observations can be made for the time step, where one even needs $\tau <\epsilon^2$ in general.
		The reason is that $u_\epsilon$ is highly oscillatory in space \emph{and} time, see Figure \ref{oszillation}. 
		\begin{figure}
			
			\begin{subfigure}{.48\textwidth}
				\centering
				\includegraphics[width=0.7\textwidth]{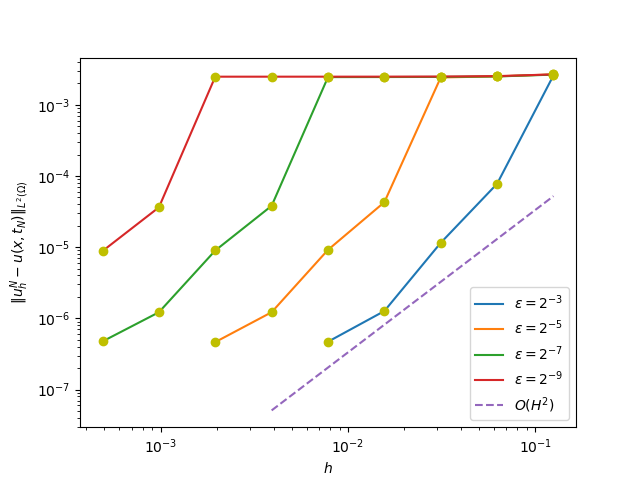}
				\caption{$L^2$-Error with respect to grid width \\ $h$ and $\epsilon$ at time $t = 1$. }
				\label{Model_L2}
			\end{subfigure}
			\begin{subfigure}{.48\textwidth}
				\centering
				\includegraphics[width=0.7\textwidth]{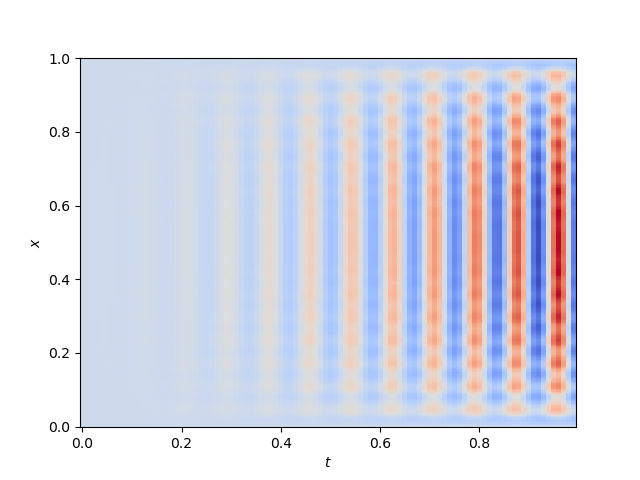}
				\caption{Illustration of $u_\epsilon$ }
				\label{oszillation}
			\end{subfigure}
	\end{figure}
	\end{Bsp}

	To tackle the outlined challenges, various multiscale methods have been proposed. Focusing on approaches for parabolic space-time multiscale problems, examples include generalized multiscale finite element methods \cite{chung2018generalized}, non-local multicontinua schemes \cite{hu2021},  high-dimensional (sparse) finite element methods \cite{tan2019high}, an approach based on an appropriate global coordinate transform \cite{owhadi2007}, a method in the spirit of the Variational Multiscale Method and the Localized Orthogonal Decomposition \cite{ljung2021} as well as optimal local subspaces \cite{schleuss2022,schleuss2022time}.
	As already mentioned, we consider locally periodic problems in space and time here. Hence, we employ the Heterogeneous Multiscale Method (HMM), first induced by E and Enquist \cite{weinan2003heterognous}, see also the reviews \cite{article,abdulle2012heterogeneous}. The HMM has been successfully applied to various time-dependent problems such as (nonlinear) parabolic problems \cite{abdulle2016numerical,abdulle2016fehmmparabolicmonotone,abdulle2015linearized,abdulle2012rungekutta}, time-dependent Maxwell equations \cite{freese2022heterogeneous,hochbruck2019,Hochbruck2017} or the heat equation for lithium ion batteries \cite{veszelkaanwendung}.
	We use the finite element version of the HMM, but note that other discretization types such as discontinuous Galerkin schemes are generally possible as well.
	
	The present contribution is inspired by \cite{ming2007analysis}, which considers the same parabolic model problem and analyzes a semi-discrete HMM for it. Precisely, the microscopic cell problems are solved analytically in \cite{ming2007analysis}. 
	Our main contribution is to propose a suitable discretization of these cell problems and to show rigorous error estimates for the resulting fully discrete HMM.
	A particular challenge for the estimate is to balance the order of the mesh size and the time step on the one hand and the period $\epsilon$ on the other hand.
	Further, we illustrate our theoretical results with numerical experiments and thereby underline the applicability of the method.
	
	The paper is organized as follows.
	In Section \ref{Multi}, we introduce the setting and present the main homogenization results.  
	In Section \ref{FEHMM}, we derive the fully discrete finite element heterogeneous multiscale method. The error of the macroscopic discretization is estimated in Section \ref{Errmacro} and the error arising from the microscopic modeling is investigated in Section \ref{eHMM}. 
	Finally, numerical results are presented in Section \ref{result}.
	
	\section{Setting \label{Multi}}
	In this section, we present our model problem and the associated homogenization results. Throughout the paper, we use standard notation on function spaces, in particular the Lebesgue space $L^2$, the Sobolev spaces $H^1$ and $H^1_0$, as well as Bochner spaces for time-dependent functions. 
	We denote the $L^2$-scalar product (w.r.t to space) by $\langle \cdot, \cdot\rangle_0$ and the $L^2$-norm by $\|\cdot\|_0$. Furthermore, we mark by $\#$ spaces of periodic functions. Let $X_{\#}(\Omega)$ be such a space for an arbitrary $\Omega \subset \mathbb{R}^d$, then the subspace $X_{\#,0}(\Omega) \subset X_{\#}(\Omega)$ consist of all functions whose integrals over $\Omega$ is $0$.
	
	\subsection{Model problem}
	Let $\Omega \subset \mathbb{R}^d$ be a bounded Lipschitz domain, $T>0$ the final time and $Y:= (-\frac{1}{2},\frac{1}{2})^d$. We consider the following parabolic problem
	\begin{align}
		\label{Pro}
		\begin{cases}
			\dfrac{\partial u^{\epsilon}}{\partial t} - \nabla \cdot (a^{\epsilon}(t,x) \nabla u^{\epsilon})  &= f(t,x)  \quad x \in \Omega, \quad t \in (0,T) \\[0.2cm]
			u^{\epsilon}(0,x) &= u_0(x) \quad \; \,  x \in \Omega \\
			u^{\epsilon}(t,x) &= 0 \quad \quad  \,\quad  \, x \in \partial\Omega, \quad t \in (0,T),
		\end{cases}
	\end{align}
	where $ f \in L^2((0,T),L^2(\Omega)) $ and $ u_0 \in L^2(\Omega)$. $a^{\epsilon}$ is the time-space  multiscale coefficient as introduced in Section \ref{Motivation} and is defined by the matrix-valued function $a(t,x,s,y) \in C([0,T] \times \bar{\Omega} \times [0,1] \times \bar{Y}, \R^{d \times d}_{sym})$.
	The function $a$ is $(0,1)\times Y$-periodic with respect to $s$ and $y$, furthermore it is coercive and uniformly bounded, in particular this means that there are constants  $\Lambda, \lambda > 0$, such that for all $\xi, \eta \in \mathbb{R}^d$: 
	\begin{align*}
		\eta \cdot a(t,x,s,y) \xi \leq \Lambda \, |\eta|\,| \xi| \quad \text{ und } \quad \xi \cdot a(t,x,s,y) \xi \geq \lambda |\xi|^2
	\end{align*}
	for all $(t,x,s,y) \in [0,T] \times \overline{\Omega} \times [0,1] \times Y$. Further, we assume that $a$ is Lipschitz continuous in $t$ and $x$.
	
	\subsection{Homogenized Problem}
	Analytical homogenization results for \eqref{Pro} were obtained in \cite{papanicolau1978asymptotic,tan2019high}. For ease of presentation, we follow the traditional approach of asymptotic expansions here, but we emphasize that the same results are obtained with the more recent approach of time-space multiscale convergence as in \cite{tan2019high}, which is a generalization of two-scale convergence.
	Based on the multiscale asymptotic expansion 
	\begin{align}
		\label{asymtotisch}
		u^{\epsilon}(t,x) = U_0(t,x,\frac{t}{\epsilon^2},\frac{x}{\epsilon}) +  \epsilon U_1(t,x,\frac{t}{\epsilon^2},\frac{x}{\epsilon}) +   \epsilon^2 U_2(t,x,\frac{t}{\epsilon^2},\frac{x}{\epsilon}) + \dots ,
	\end{align}
	it is shown that $U_0$ solves the \textit{homogenized problem}
	\begin{align}
		\label{homo}
		\begin{cases}
			\dfrac{\partial U_0}{\partial t} - \nabla \cdot (A_{0} \nabla U_0)  &= f \quad \text{ in 		} (0,T) \times  \Omega  \\
			U_0 &= 0 \quad \text{ on } (0,T) \times  \partial\Omega \\
			U_0(0,\cdot) &= u_0 \quad  \text{ in } \Omega.  \\
		\end{cases}
	\end{align}
	Here, the homogenized coefficient $A_0$ is defined by
	\begin{align}
		\label{A_0_Ursprung}
		A_0^{ij}(t,x) = \int_{0}^{1}\int_{Y} \sum_{k = 1}^d a_{ik}(t,x,s,y)\bigl(\delta_{jk} +				\dfrac{\partial \chi^{j}}{\partial y_k}(t,x,s,y)\bigr) dyds,
	\end{align}
	where $\delta_{jk}$ denotes the Kronecker delta.
	The function $\chi^i \in L^2((0,T) \times \Omega \times (0,1), H^1_{\#,0}(Y)) \cap L^2((0,T) \times \Omega ,H^1_{\#}((0,1),H^{-1}_{\#,0}(Y)))$  solves the cell problem
	\begin{align}
		\label{cell}
		\begin{cases}
			\dfrac{\partial \chi^{i}}{\partial s} - \,\nabla_{y} \cdot(a(e^{i} + \nabla_{y}				\chi^{i})) 	= 0  \quad \text{ in }  (0,1) \times Y,\\
			\chi^{i}(t,x,s, \cdot) \quad Y\text{-periodic for all } t,x,s, \\
			\chi^{i}(t,x, \cdot,y) \quad (0,1)\text{-periodic for all } t,x,y.
		\end{cases}
	\end{align}
	Using these $\chi^i$, $U_1$ in the asymptotic expansion \eqref{asymtotisch} can be written as 
	\begin{align*}
		U_1(t,x,s,y) = \sum_{i = 1}^{d} \frac{\partial U_0}{\partial x_i}(t,x) \chi^{i}(t,x, 				s,y).
	\end{align*}
	\cite[Chapter 2, Section 1.7]{papanicolau1978asymptotic}   shows in Theorem 2.1 and Theorem 2.3 that
	\begin{align*}
		|\!| u^{\epsilon} - U_0 - \epsilon U_1 |\!|_{L^2((0,T),H^1_0(\Omega))} \rightarrow 0  \text{ for } \epsilon \rightarrow 0 .
	\end{align*}
	We call $U_0$ the \textit{homogenized solution}. 
	$U_0$ describes the macroscopic behavior of $u^{\epsilon}$, because $U_0$ only depends on the macroscopic scale $x$. $U_1$ is called \textit{the first-order corrector}. 
	\begin{Bemerkung}
		$A_0$ is not symmetric in $\R^d$ for $d>1$ in general since 
		\begin{multline}
			\label{A0symm}
			A_0^{ij}(t,x) = \int_{0}^{1}\int_{Y} \sum_{l = 1}^d \sum_{k = 1}^d \delta_{il} a_{lk}   		(t,x,s,y)(\delta_{jk} +\dfrac{\partial \chi^{j}}{\partial y_k}(t,x,s,y)) dyds \\
			= \int_{0}^{1}\int_{Y} \sum_{l = 1}^d \sum_{k = 1}^d (\delta_{il} +\dfrac{\partial 			\chi^{i}}{\partial y_l}(t,x,s,y)) a_{lk}(t,x,s,y) (\delta_{jk} +\dfrac{\partial 				\chi^{j}}{\partial y_k}(t,x,s,y)) dyds  \\ +  \int_{0}^{1}\int_{Y} \chi^i(t,x,s,y) \dfrac{\partial 			\chi^j}{\partial s} (t,x,s,y) dyds.
		\end{multline}
		The last term does not vanish in general, but it is zero for $i = j$ due to integration by parts and the time-periodicity of $\chi^i$.
	\end{Bemerkung}
	In the following we reformulate $A_0$ in a way which we use to derive the discretized problem later. We transform the reference cell $(0,1) \times Y$ to a general cell $(t,t + \epsilon^2) \times \{x_0\} + I_{\epsilon}$ with $I_\epsilon\coloneqq \epsilon Y$ for $t \in [0,T)$ and $ x_0 \in \Omega$ fixed. Application of the chain and transformation rule allows us to write
	\begin{align}
		\label{A_0_finial}
		A_0(t,x ) &= \int_{0}^{1}\int_{Y}   a(t,x,s,y) 					(\operatorname{Id}_d +  D_y \chi(t,s,y)) dyds \\
		&=  \dfrac{1}{\epsilon^2|I_{\epsilon}|}\int_{t}^{t+ \epsilon^2}\int_{\{x_0\} + 						I_{\epsilon}} 	a\Big(t,x,\dfrac{s}{\epsilon^2},\dfrac{y}{\epsilon}\Big)\Big(\operatorname{Id}_d +  D_y \chi\Big(t,x,\dfrac{s}			{\epsilon^2},\dfrac{y}{\epsilon}\Big)\Big)\,  dyds.
	\end{align}
	
	\section{The finite-element heterogeneous multiscale method (FE-HMM) \label{FEHMM}}
	Based on the results of Section \ref{Multi}, we want to compute an approximation of the homogenized solution $U_0$ based on the Finite-Element Heterogeneous Multiscale Method (FE-HMM).
	In \cite{ming2007analysis}, this method was already introduced, but it was assumed that the cell problems \eqref{cell} could be solved exactly/analytically. The main aim of this section is to introduce also  the (microscopic) discretization of the cell problems, allowing for a fully discrete method.
	Further, we also account for the non-symmetry of $A_0$. This leads to a slightly different formulation in comparison to \cite{ming2007analysis} where the symmetric part of $A_0$ was considered throughout.
	In the following, we will derive the full method step by step, which is on the one hand hopefully instructive for the readers to understand the final formulation and on the other makes it easier to follow the error estimates in the following sections.
	
	We start with the discretized macro problem. For the spatial discretization we use linear finite elements based on a triangulation $\mathcal{T}_H$ and for the time discretization we use the implicit Euler method.
	Precisely, let $V_{H} \subset H_0^1(\Omega)$ be the space of all piecewise linear functions which are zero on $\partial \Omega$  and let $\tau = T/N$ be the time step size. For $1 \leq n\leq N$ we set $ t_n = n \tau$. Further, we define $U_H^0 := Q_Hu_0$, where $Q_H:L^2(\Omega)\to V_H$ is the $L^2$-projection. 
	 
	Let $U_H^n$ then be the solution of the discretized equation
	\begin{align}
		\label{diskret1}
		\Bigl\langle\dfrac{\overline{\partial} U_H^n}{\partial t}, \Phi_H \Bigr\rangle_{0} + B[t_n,U_H^n, \Phi_H] = \bigl\langle f^n,\Phi_H\bigr\rangle_{0} \, \text{ for all } \Phi_H \in V_{H},
	\end{align}
	where $f^n(x)= f(t_n,x)$ and $ \dfrac{\bar{\partial} U_H^n}{\partial t} = (U^n_H - U^{n-1} _H)/ \tau$. 
	Here, the discrete bilinear form $B[t_n,\cdot,\cdot]$ is defined for any $\Phi_H, \Psi_H \in V_H$ via
	\begin{align}	
		B[t_n,\Phi_H,\Psi_H] &:= \int_{\Omega} \nabla \Psi_H(x)\cdot A_0(t_n,x) \nabla \Phi_H(x) 		dx \notag \\ &= \sum_{K \in \mathcal{T}_H} \int_{K} \nabla \Psi_H(x)\cdot A_0(t_n,x) \nabla 			\Phi_H(x) dx \notag \\ &\approx \sum_{K \in \mathcal{T}_H} |K| \nabla \Psi_H(x_K)\cdot \A0nK 	\nabla 	\Phi_H(x_K).
		\label{quadratur_1}
	\end{align}
	In the last step we approximated the integral with a quadrature formula, where $x_K$ denotes hte barycenter of $K\in \mathcal T_H$. 
	If we now consider the individual summands, we could calculate $\A0nK$ starting from equation (\ref{A_0_finial}). However, this would have several disadvantages. First, we would have to compute $A_0(t_n,x_K) $ for all time points $t_n$, which would require a lot of memory depending on the time step size. Furthermore, we want to change the boundary conditions later, which is not possible with this approach. Therefore, the idea is to compute $\nabla\Psi_H(x_K)\cdot \A0nK \nabla \Phi_H(x_K)$ directly. For this, set $I_{\epsilon,K} := \{x_K\} + I_{\epsilon}$ and use reformulation \eqref{A_0_finial} to give
	\begin{align}
		&\nabla \Psi_H(x_K) \cdot \A0nK  \nabla \Phi_H(x_K) \notag \\  &=  \dfrac{1}{\epsilon^2|				I_{\epsilon}|}\int_{t_n}^{t_n+ \epsilon^2} \int_{ I_{\epsilon,K}}
		\nabla\Psi_{H \mid_{ I_{\epsilon,K}}}(x_K)  \cdot \overbrace{a(t_n,x_K,\dfrac{t}{\epsilon^2},\dfrac{x}{\epsilon})}^{ : = \a0n (t,x) }\notag  \\  & \qquad \qquad \qquad \qquad \qquad \qquad \qquad \qquad \qquad \nabla(\Phi_{H 					\mid_{I_{\epsilon,K}}}(x_K) + \underbrace{\nabla\Phi_{H \mid_{I_{\epsilon,K}}}(x_K)\epsilon 				\chi(t_n,x_K\dfrac{t}{\epsilon^2},\dfrac{x}{\epsilon})}_{=: \tilde{\Phi}})dxdt \notag \\ &= 			\dfrac{1}{\epsilon^2|I_{\epsilon}|}\int_{t_n}^{t_n+ \epsilon^2} \int_{ I_{\epsilon,K}} 	 \nabla\Psi_{H \mid_{ I_{\epsilon,K}}}(x_K) \cdot 		 \a0n (t,x)  \nabla \phi_{\#}^{\epsilon}			(t_n,x_K,\dfrac{t}			{\epsilon^2},\dfrac{x}{\epsilon}) dx dt,
		\label{Substitution_Phi_Psi}
	\end{align}	
	where
	\begin{align}
		\label{space_phi}
		\phi^{\epsilon}_{\#} = \Phi_{H \mid_{I_{\epsilon,K}}} +   \tilde{\Phi} \in V_{H} + 	X((t_n,t_n + \epsilon^2),I_{\epsilon,K}), 
	\end{align}
	with
	\begin{align*}
		X((t_n,t_n + \epsilon^2),I_{\epsilon,K}) := L^2((t_n,t_n + \epsilon^2),H^1_{\#,0}(I_{\epsilon,K} )) \cap H^1_{\#}((t_n,t_n + \epsilon^2),H^{-1}_{\#,0}(I_{\epsilon,K} )).
	\end{align*}
	$\phi^{\epsilon}_{\#}$ solves the equivalent cell problem
	\begin{align}
		\label{cell_ae}
		\begin{cases}
			\dfrac{\partial \phi^{\epsilon}_{\#}}{\partial t} - \,\nabla	\cdot(							\a0n  \nabla	\phi_{\#}^{\epsilon} )= 0  \quad &\text{in }  (t_n,t_n + 		\epsilon^2) \times 	I_{\epsilon,K} \\
			\phi_{\#}^{\epsilon}(t,x,s, \cdot)-\Phi_{H \mid_{I_{\epsilon,K}}} \qquad &\text{periodic on } \partial I_{\epsilon,K}	\\
			\phi_{\#}^{\epsilon}(t,x, \cdot,y)-\Phi_{H \mid_{I_{\epsilon,K}}} \qquad  &\text{periodic on } \partial(t_n, t_n+\epsilon^2).
		\end{cases}
	\end{align}
	We can thus give the first discretization for the bilinear form in \eqref{diskret1}
	\begin{align*}
		B_{H,\#}[t_n,\Phi_H,\Psi_H]  
		&:= \sum_{K \in \mathcal{T}_H} \dfrac{|K|}{|I_{\epsilon}|\epsilon^2} 				\int_{t_n}^{t_n+ \epsilon^2} \int_{ I_{\epsilon,K}} 	\nabla\PhiH(x_K)\cdot \a0n(t,y)	\nabla							\phi^{\epsilon}_{\#}(t,y)dydt \\
		&\approx 
		\sum_{K \in \mathcal{T}_H} \int_K \dfrac{1}{|I_{\epsilon}|\epsilon^2} 				\int_{t_n}^{t_n+ \epsilon^2} \int_{ I_{\epsilon,K}} 	\nabla\PhiH(x)\cdot \a0n(t,y)	\nabla							\phi^{\epsilon}_{\#}(t,y)dydtdx 
	\end{align*}
	Note that $ B_{H,\#}$ is generally not symmetric.
	
	In practice, the period may be known only approximately. Therefore we consider the case with cell side length $\delta > \epsilon$ and cell time $\sigma > \epsilon^2$ and where the two terms $\frac{\sigma}{\epsilon^2} $,
	$\frac{\delta}{\epsilon}$ are not integers. Thus, the periodic boundary conditions no longer hold (see \cite[p.164]{article} for the stationary case). In this case, we need to find alternative boundary and initial values.
	We approximate $\nabla \Psi_H \cdot A_0(t_n) \nabla \Phi_H$ by replacing $\epsilon$ by $\delta$ and $\epsilon ^2$ by $\sigma$ in \eqref{Substitution_Phi_Psi}, and also $X$ by 
	\begin{align*}
		L^2((t_n,t_n + \sigma),H^1_{0}(I_{\delta,K} )) \cap H^1((t_n,t_n + \sigma),H^{-1}_{0}(I_{\delta,K} ))
	\end{align*}
	in \eqref{space_phi}. 
	Then we approximate
	\begin{align*}
		\nabla \Psi_H(x_K) \cdot A_0(t_n,x_K) \nabla \Phi_H(x_K) &\approx  \dfrac{1}{\sigma|I_{\delta}|}		\int_{t_n}^{t_n+ \sigma}\int_{I_{\delta,K}} (\nabla\Phi_{H }(x_K)\cdot \a0n			(t,x)\nabla\phi^{\epsilon}(t,x))dxdt \\
		&:=\dfrac{1}{\sigma|I_{\delta}|}		\int_{t_n}^{t_n+ \sigma}\int_{I_{\delta,K}} (\nabla\Phi_{H \mid_{ I_{\delta,K}}}(x_K)\cdot \a0n			(t,x)\nabla\phi^{\epsilon}(t,x))dxdt,
	\end{align*}
	where $\phi^{\epsilon}$ solves the  initial value problem
	\begin{align}
		\label{cell_anfangs}
		\begin{cases}
			\dfrac{\partial \phi^{\epsilon}}{\partial t} - \,\nabla	\cdot(									 \a0n \nabla \phi^{\epsilon} ) &= 0  \quad \text{ in }  (t_n,t_n + \sigma) 			\times 	I_{\delta,K} \\
			\phi^{\epsilon}  &= \Phi_H \quad \text{ on }  (t_n, t_n + \sigma) \times	\partial 			I_{\delta,K}\\
			\phi^{\epsilon}_{\mid_{t = t_n} } &= \Phi_H.
		\end{cases}
	\end{align}
	This means that we replace periodic boundary conditions by Dirichlet  ones and the time ``boundary value problem'' by an initial value problem.
	
	In the following, we consider the bilinear form resulting from the above approximation. We set
	$ \mathcal{Q}_{n,K}:= (t_n,t_n + \sigma) 	\times 	I_{\delta,K}$ and define	
	\begin{align}	
		B_H[t_n,\Phi_H,\Psi_H]  &:=
		\sum_{K \in \mathcal{T}_H} |K| 				\nabla \Psi_H(x_K) \cdot A_H(t_n,x_K) \nabla \Phi_H(x_K),\label{quadratur} \\
		&=  \sum_{K \in \mathcal{T}_H} \int_K				\nabla \Psi_H(x) \cdot A_H(t_n,x_K) \nabla \Phi_H(x) dx, \label{B_H} 
	\end{align}
	where
	\begin{align*}
		\nabla \Psi_H \cdot A_H(t_n, x_K) \nabla \Phi_H 
		 &:=\dfrac{1}{|\mathcal{Q}_{n,K}|}  \int_{ \mathcal{Q}_{n,K}} \nabla	 \Psi_H(x_K)\cdot \a0n(t,x)		\nabla					\phi^{\epsilon}(t,x)dxdt.
	\end{align*}
	To finally get the fully discrete method we consider a triangulation $T_{\tilde{h}}$ of the unit cell $Y$ and the resulting triangulation $T_h(I_{\delta,K})$ of the shifted cell with the finite element space $\Vmicro \subset H_0^1(I_{\delta,K})$, which consist of all piecewise polynomials of order $p \geq 2$. We stress that the mesh size $h$ is meant with respect to the scaled triangulation  $T_h(I_{\delta,K})$.
	Let $1\leq k \leq N_{cell}$, $\theta = \frac{\sigma}{N_{cell}}$ and $s_k = k \theta$. For any $\Phi_H \in V_H$ seek $\phidisc \in \Phi_H + \Vmicro$ as the unique solution of the discrete cell problem
	\begin{align*}	
		\int_{I_{\delta,K}} \dfrac{\overline{\partial} \phidisc}{\partial t}z_h  + \nabla z_h \cdot \a0n  \nabla \phidisc dx = 0 \quad \text{for all }z_h \in \Vmicro
	\end{align*}
	 with $\phi_{h,0}^{\epsilon} = \Phi_H$.
	
	We consider the following  bilinear form, which we get from the approximations above	
	\begin{align}	
		\Bhh[t_n,\Phi_H,\Psi_H] & :=
		\sum_{K \in \mathcal{T}_H} |K| 				\nabla \Psi_H(x_K) \cdot \Ahh(t_n,x_K) \nabla \Phi_H(x_K)\label{quadratur_disc}  \\
		&\; =  \sum_{K \in \mathcal{T}_H} \int_K				\nabla \Psi_H(x) \cdot \Ahh(t_n,x_K) \nabla \Phi_H(x) dx, \label{B_H_disc} 
	\end{align}
	where
	\begin{align*}
		\nabla \Psi_H \cdot \Ahh(t_n,x_K) \nabla \Phi_H &:=
		\dfrac{1}{\sigma|I_{\delta}|} \Bigl(\dfrac{\theta}{2} \int_{ I_{\delta,K}} \nabla \Psi_H(x_K)\cdot \a0n(t_n,x) \nabla \phi^{\epsilon}_{h,0} dx \\ &\quad +  \theta \sum_{ k= 1}^{N_{cell}-1}\int_{ I_{\delta,K}} \nabla \Psi_H(x_K) \cdot \a0n(t_n + s_k,x) \nabla \phi^{\epsilon}_{h,k}dx   \\ & \quad + \dfrac{\theta}{2} \int_{ I_{\delta,K}} \nabla  \Psi_H(x_K)\cdot \a0n(t_n + s_{N_{cell}},x) \nabla\phi_{h,n}^{\epsilon} dx \Bigr).
	\end{align*}
	\begin{Bemerkung}
		In the above definition of $A_{H,h}$ we used the trapezoidal rule for the  approximation for the time integral, which is consistent with our numerical experiments below. 
		We emphasize that the choice of other quadrature rules is equally possible. In practice, also the spatial integral over $I_{\delta, K}$ is approximated by a quadrature rule.
	\end{Bemerkung}
	We reformulate the discrete homogeneous equation by substituting $B_{H,h}[t_n,\Phi_H,\Psi_H]$ into (\ref{diskret1}): Seek $U_H^n\in V_H$ such that 
	\begin{align}
		\label{diskret2}
		\Bigl\langle\dfrac{\bar{\partial} U_H^n}{\partial t}, \Phi_H \Bigr\rangle_0 + 					\Bhh[t_n,U_H^n,		\Phi_H] = \bigl\langle f^n,\Phi_H\bigr\rangle_0  \text{ for all		} \Phi_H \in V_{H},
	\end{align}
	where again $U_H^0=Q_H u_0$.
	\begin{Bemerkung} 
		The restriction to use only 	piecewise linear functions for the macrodiscretization  is
		important for proving  the claimed error bounds later. For finite elements of higher degree, one approach is to consider the linearization $\Phi_{H,lin}$ of finite element functions $\Phi_H \in V_H$, where 
		\begin{align*}
			\Phi_{H,lin} : = \sum_{K \in \mathcal{T}_H} \Phi_{H,lin,K} = \sum_{K \in \mathcal{T}_H} \Phi_{H}(x_K) + (x -x_K) \cdot\nabla \Phi_H(x_K).
		\end{align*}
		We refer to \cite{article} for details in the stationary case.
	\end{Bemerkung}

In the following, we will show well--posedness as well as error estimates for the FE-HMM. Note that for the well--posedness, it is sufficient to show stability for the macrodiscretization since the microproblem is ``only'' used to calculate the homogenized coefficient $A_{H,h}$ or the discrete bilinear form $B_{H,h}$, respectively. 
	
	\section{Error estimation of the macrodiscretization \label{Errmacro}}
	We start with a stability result which we will use to show coercivity and boundedness of $\Bhh$. The statement and the proof are similar to \cite[Lemma 2.1]{ming2007analysis}. However, we consider here the time-discretized case.
	\begin{Lemma}
		\label{4stab}
		Let $\Omega \subset \mathbb{R}^d$ be a bounded domain, $ T >0$ and $\Phi$ a linear function. Further, let $\varphi$ be a solution to the following problem.
		\begin{align}
			\label{Diff_Allgemein}
			\begin{cases}
				\dfrac{\partial \varphi}{\partial t} - \,\nabla	\cdot( a(t,x) \cdot  \nabla	\varphi ) &= 0  \quad \text{ in }  (0,T] \times 	\Omega \\
				\varphi &= \Phi \; \; \; \text{ on }  (0,T] \times 	\partial\Omega 	\\
				\varphi\!\mid_{ t = 0} &= \Phi,
			\end{cases}
		\end{align}
		where $a(t,x) = (a_{ij}(t,x))_{i,j = 1\dots d}$ fulfills the following conditions
		\begin{align*}
			\lambda \operatorname{Id}_d  \leq a(t,x) \leq \Lambda \operatorname{Id}_d  \quad \text{  a.e. on } (0,T] \times \Omega,
		\end{align*}
		where $\operatorname{Id}_d: \mathbb{R}^d \to \R^d$ is the $d$-dimensional unit matrix.
		Let $V_{h}^p \subset H_0^1(\Omega)$ be the space of all piecewise polynomials of degree $p$ on $\Omega$ using the simplicial mesh $\mathcal T_h$.
		For $1 \leq n \leq N$ we define $\theta = T/N$, $t_n = n \theta $ and let $\varphi_h^n \in \Phi + V_h^p$ be the weak solution of the discrete problem
		\begin{align*}	
			\int_{\Omega} \dfrac{ \varphi^n_h - \varphi^{n-1}_h}{\theta}z_h + \nabla z_h \cdot a(t_n,x)  \nabla\varphi^n _h(x)  dx = 0
		\end{align*}
		for all $z_h \in V_h$ and  $\varphi_h^0 = \Phi$. 
		Then it holds for all $ n \in \{0,\dots ,N\}$:
		\begin{align}
			\label{Ungleich_Allgemein1}
			|\!| \nabla \Phi|\!|_{0} &\leq |\!| \nabla \varphi_h^n |\!|_{0} \\ \label{Ungleich_Allgemein2}C T^{1/2}|\!| \nabla \Phi|\!|_{0} &\geq \Bigl(\theta \sum_{i=0}^{N}|\!| \nabla( \varphi_h^i- \Phi) |\!|^2_{0} \Bigr)^{1/2}.
		\end{align}
	\end{Lemma}	
	\begin{proof}
		Let $ n \in \{0,\dots ,N \}$ be arbitrary. Since $ \varphi_h^n = \Phi$ on $\partial \Omega$ and $ \nabla \Phi$ is constant, we get with partial integration 
		\begin{align*}
			\int_{\Omega}\nabla\bigl(\varphi_h^n(x) - \Phi(x)\bigr) \cdot \nabla \Phi(x)dx	&= -\int_{\Omega} \bigl(\varphi_h^n(x) - \Phi(x)\bigr) \Delta \Phi(x) dx  \\& \quad  + \oint_{\partial\Omega} \bigl(\varphi_h^n(x) -\Phi(x)\bigr) \nabla \Phi(x)\cdot\nu d\nu \\  &= 0 
		\end{align*}
		and, hence,
		\begin{align*}
			\int_{\Omega} |\nabla \varphi_h^n(x)|^2dx = \int_{\Omega} |\nabla \Phi(x)|^2dx + \int_{\Omega} |\nabla\bigl(\varphi_h^n(x) - \Phi(x)\bigr)|^2 dx.
		\end{align*}
		The first inequality (\ref{Ungleich_Allgemein1}) is thus proved.
		For the second inequality, we use that $\varphi$ is the weak solution of (\ref{Diff_Allgemein}). We choose $\varphi_h^ n - \Phi $ as the test function to obtain
		\begin{multline*}
			\sum_{n = 1}^N \Big[	\int_{\Omega} (\varphi_h^n(x) - \varphi_h^{n-1}(x)) (\varphi_h^ n(x) - \Phi(x)) dx + \;\theta\int_{\Omega}  \nabla (\varphi_h^ n(x) - \Phi) \cdot a(t_n,x) \nabla (\varphi_h^ n(x) - \Phi (x)) dx \Big] \\* = \sum_{n = 1}^N \theta \int_{\Omega}  \nabla (\varphi_h^ n(x) - \Phi(x)) \cdot a(t_n,x) \nabla \Phi(x)  dx .
		\end{multline*}
		From the Cauchy-Schwarz inequality and boundedness of $a$ it follows
		\begin{align*}
			\sum_{n = 1}^N  \int_{\Omega}&  \nabla (\varphi_h^ n(x) - \Phi (x)) \cdot a(t_n,x) \nabla \Phi(x)  dx \\ &\leq   \Bigl( \sum_{n = 1}^N \int_{\Omega} \nabla (\varphi_h^ n(x) - \Phi(x) )  \cdot a(t_n,x) \nabla (\varphi_h^ n(x) - \Phi (x)) dx \Bigr)^{1/2} \\ & \qquad \qquad \quad \qquad \Bigl(\sum_{n = 1}^N \int_{\Omega}  \nabla \Phi(x) \cdot a(t_n,x) \nabla \Phi(x) dx \Bigr)^{1/2} \\
			&\leq \Lambda^{1/2} \sqrt{N} \| \nabla \Phi\|_0   \Bigl( \sum_{n = 1}^N\int_{\Omega} \nabla (\varphi_h^ n(x) - \Phi(x) )  \cdot a(t_n,x) \nabla (\varphi_h^ n(x) - \Phi (x)) dx \Bigr)^{1/2}.
		\end{align*}
		Inserting this inequality into the equation above and using that 
		\begin{align*}
			\sum_{n= 1}^N\int_{\Omega} (\varphi_h^n(x) - \varphi_h^{n-1}(x)) (\varphi_h^ n(x) - \Phi(x)) dx \geq \frac{1}{2}(\|\varphi_h^N - \Phi\|_0^2 - \|\varphi_h^0 - \Phi\|_0^2) \geq 0,
		\end{align*}
		we finally obtain
		\begin{align*}
			\sum_{n = 1}^N \theta &\int_{\Omega} \nabla (\varphi_h^ n(x)  - \Phi(x)  ) \cdot a(t_n,x) \nabla (\varphi_h^ n(x)  - \Phi (x)) dx  \\
			&\leq C \theta^{1/2} \| \nabla \Phi\|_0  \Bigl(\sum_{n = 1}^N \theta \int_{\Omega} \nabla (\varphi_h^ n(x)  - \Phi(x)  ) \cdot a(t_n,x) \nabla (\varphi_h^ n(x)  - \Phi (x)) dx \Bigr)^{1/2}
		\end{align*}
		Finally dividing by $\bigl(\sum_{n = 1}^N \theta \int_{\Omega} \nabla (\varphi_h^ n - \Phi )(x)  \cdot a(t_n,x) \nabla (\varphi_h^ n - \Phi (x)) dx \bigr)^{1/2}$ finishes the proof.
	\end{proof}
	Using this lemma, we now show the boundedness and coercivity of the discretized bilinear form $B_{H,h}$.
	\begin{Lemma}
		For all $ n \in \{1,\dots, N\}$, $\Bhh[t_n,\cdot,\cdot]: V_H\times V_H \rightarrow \mathbb{R}$ is a coercive and bounded bilinear form.
	\end{Lemma}
	\begin{proof}
		Using the definition of  $\phi^{\epsilon}_{h,k}$, we obtain with Lemma \ref{4stab} and $\sigma=N_{cell} \theta$
		\begin{align*}
			&\!\!\!\!\nabla \Psi_H \cdot \Ahh(t_n,x_K) \nabla \Phi_H \\ &= \dfrac{1}{\sigma|I_{\delta}|} \Bigl(\dfrac{\theta}{2} \int_{I_{\delta,K}} \nabla \Psi_H(x_K)\cdot \a0n(t_n,x) \nabla \phi^{\epsilon}_{h,0} dx \\ &\qquad   \qquad +  \theta \sum_{ k= 1}^{N_{cell}-1}\int_{I_{\delta,K}} \nabla \Psi_H(x_K) \cdot \a0n(t_n + s_k,x) \nabla \phi^{\epsilon}_{h,k}dx   \\ & \qquad  \qquad + \dfrac{\theta}{2} \int_{ I_{\delta,K}} \nabla  \Psi_H(x_K)\cdot \a0n(t_n + s_{N_{cell}},x) \nabla\phi_{h,n}^{\epsilon} dx \Bigr) \\
			&\leq \dfrac{1}{\sigma|I_{\delta}|} C\Bigl(\dfrac{\theta}{2} |\!|\nabla \Psi_H|\!|_{L^2(I_{\delta,K})} \bigl(|\!| \nabla( \phi^{\epsilon}_{h,0}  - \Phi_H)|\!|_{L^2(I_{\delta,K})} + \! |\!|\nabla \Phi_H|\!|_{L^2(I_{\delta,K})}\bigr)  \\* & \qquad  \qquad \quad +   \theta \sum_{ k= 1}^{N_{cell}-1} |\!| \nabla \Psi_H|\!|_{L^2(I_{\delta,K})} \bigl(|\!| \nabla(\phi^{\epsilon}_{h,k}- \Phi_H)|\!|_{L^2(I_{\delta,K})}   + |\!|\nabla \Phi_H|\!|_{L^2(I_{\delta,K})}\bigr) \\* & \qquad  \qquad \quad + \dfrac{\theta}{2}  |\!| \nabla  \Psi_H |\!|_{L^2(I_{\delta,K})} \bigl(|\!|\nabla(\phi_{h,n}^{\epsilon} - \Phi_H)|\!|_{L^2(I_{\delta,K})}   + |\!|\nabla \Phi_H|\!|_{L^2(I_{\delta,K})}\bigr) \Bigr)
			\\ &\leq \dfrac{1}{\sigma|I_{\delta}|} C   \Big[ \Bigl(\theta\sum_{ k= 0}^{N_{cell}} |\!| \nabla \Psi_H|\!|^2_{L^2(I_{\delta,K})}  \Bigr)^{1/2}  \Bigl(\theta\sum_{ k= 0}^{N_{cell}}  |\!| \nabla(\phi^{\epsilon}_{h,k}- \Phi_H)|\!|^2_{L^2(I_{\delta,K})} \Bigr)^{1/2} \\ & \qquad  \qquad \qquad  + \theta\sum_{ k= 0}^{N_{cell}} |\!|\nabla \Psi_H|\!|_{L^2(I_{\delta,K})} |\!|\nabla \Phi_H|\!|_{L^2(I_{\delta,K})}
			\Big]
			\\ &\leq C\frac{1}{|I_\delta|}\|\nabla \Psi_H\|_{L^2(I_{\delta,K})}\|\nabla \Phi_H\|_{L^2(I_{\delta,K})}=C |\nabla \Psi_H(x_K)||\nabla \Phi_H(x_K)|,
		\end{align*}
	where we used in the last step that $\nabla \Phi_H$ is constant.
	Summation over $K$ shows the boundedness of $B_{H,h}$.

		It remains to show the coercivity. Since $ \frac{\Ahh(t_n,x_K) - \Ahh(t_n,x_K)^T}{2}  $ is skew-symmetric, it follows that
		\begin{align*}
			\nabla \Phi \cdot \Bigl(	\dfrac{\Ahh(t_n,x_K) - \Ahh(t_n,x_K)^T}{2} \Bigr) \nabla \Phi = 0.
		\end{align*}
		For the symmetric part, we obtain with \eqref{A0symm}
		\begin{align*}
			\dfrac{A_0(t,x) + A_0(t,x)^T}{2} \! = \! \int_{0}^{1}\int_{Y} \sum_{l = 1}^d \sum_{k = 1}^d (\delta_{il} +\dfrac{\partial 			\chi^{i}}{\partial y_l}(t,x,s,y)) a_{lk}(t,x,s,y) (\delta_{jk} +\dfrac{\partial 				\chi^{j}}{\partial y_k}(t,x,s,y)) dyds.
		\end{align*}
		With the lower bound on $a_{n,K}^\varepsilon$ and Lemma \ref{4stab}, we calculate
		\begin{align*}
			\nabla \Phi_H \cdot \Ahh(t_n,x_K) \nabla \Phi_H &= 
			\nabla \Phi_H \cdot \Bigl( \dfrac{\Ahh(t_n,x_K) + \Ahh(t_n,x_K)^T}{2}  \Bigr) \nabla \Phi_H \\
			&= \dfrac{1}{\sigma|I_{\delta}|} \Bigl(\dfrac{\theta}{2} \int_{I_{\delta,K}} \! 
			\nabla \phi^{\epsilon}_{h,0} \cdot \a0n(t_n,x) \nabla \phi^{\epsilon}_{h,0} dx \\ &\qquad  \qquad +  \theta \sum_{ k= 1}^{N_{cell}-1}\int_{ I_{\delta,K}} \nabla \phi^{\epsilon}_{h,k}  \cdot \a0n(t_n + s_k,x) \nabla \phi^{\epsilon}_{h,k}dx   \\ & \qquad  \qquad  + \dfrac{\theta}{2} \int_{ I_{\delta,K}} \nabla  \phi^{\epsilon}_{h,n} \cdot \a0n(t_n + s_{N_{cell}},x) \nabla\phi_{h,n}^{\epsilon} dx \Bigr) \\
			&\geq C _{0}|\nabla \Phi_H(x_K)|^2
		\end{align*}
	and the coercivity of $B_{H,h}$ follows by summation over $K$ and the fact that $\Phi_H$ is piecewise linear.			
	\end{proof}
	The macrodiscretization is thus a usual finite element discretization with implicit Euler time stepping of a coercive and bounded parabolic (discrete) problem, which directly implies its stability.
	
	We will now present the main error estimate which is of similar form as in \cite{ming2007analysis}.
	We define the error arising from the estimation of microscopic data as 
	\begin{align*}
		e(\text{HMM}) = \max_{1 \leq k \leq n} e_k(\text{HMM}),
	\end{align*}
	where
	\begin{align*}
		e_k(\text{HMM})   &:= \max_{K \in \mathcal{T}_H}|\!| (A_0 - \Ahh)(t_k,x_K) |\!| 
		 \\* &: = \max_{K \in \mathcal{T}_H} \Biggl[ \;\max_{\substack{\Phi_H,\Psi_H \in V_H \\ |\nabla\Phi_{H|_{K}}|,|\nabla\Psi_{H|_K}| = 1 }} |\nabla\Psi_H \cdot (A_0 - \Ahh)(t_k,x_K)\nabla \Phi_H| \! \Biggr].
	\end{align*}
Note that the definition is analogous to \cite{ming2007analysis}, but includes the microscopic discretization by comparing $A_0$ with $\Ahh$ and not $A_H$.
Using the same perturbation argument as \cite{ming2007analysis} one directly obtains the main error estimate.
	\begin{Satz}
		\label{Fehler_Satz}
		Let $U_0$ and $U^n_H$ be solutions of (\ref{homo}) and 		(\ref{diskret2}), respectively. If $a$ and $U_0$ are sufficiently regular, there exists a constant $C$ independent of $\epsilon, \delta, \sigma, H , \tau $ such that
		\begin{align}
			\label{Fehler3}
			|\!| U^n_H - U_0(x,t_n)|\!|_0 + |\!|\!| U^n_H - U_0(x,t_n)|\!|\!| &\leq C( \tau +H^2 				+e(\text{HMM})),\\
			\label{Fehler4}
			|\!| U^n_H - U_0(x,t_n)|\!|_{H^1(\Omega)}  &\leq C( \tau +H+e(\text{HMM})\tau^{-\frac{1}{2}}),
		\end{align}	
	where $|\!|\!|.|\!|\!|$ is defined as
	\begin{align*}
		|\!|\!|\Phi|\!|\!| : = \Bigl(\sum_{k =1}^n \tau |\!|\nabla \Phi^k|\!|_0^2\Bigr)^{1/2}
	\end{align*}
	for all $\Phi = \{\Phi^k\}_{k=1}^n$ with $\Phi^k 			\in H_0^1(\Omega)$. 
	\end{Satz}
	
	\section{Estimation of  \textit{e}(HMM)\label{eHMM}}
	In this section, we prove the error bound for $e($HMM$)$. In contrast to \cite{ming2007analysis}, we also consider the error of the microscopic discretization, which requires additional effort. Further slight differences to \cite{ming2007analysis} arise from the lacking symmetry of $A_0$ and $A_{H,h}$.
	Instead of  estimating $e($HMM$)$ directly, we introduce  auxiliary matrices $\tilde{A}$ and $\tilde{A}_H$ and calculate the error with respect to them.	Let $\Phi_H, \Psi_H \in V_H$.
	$\tilde{A}$ is defined via
	\begin{align*}
		\nabla \Psi_H \cdot \tilde{A}(t_n,x_K)\nabla \Phi_H := \avint_{\mathcal{Q}_{n,K}} \nabla \Psi_H(x) \cdot \a0n(t,x) \nabla \phi^{\epsilon}_{\#}(t,x) dxdt ,
	\end{align*}
	where $\phi^{\epsilon}_{\#}$ solves \eqref{Substitution_Phi_Psi}.
	Note that $\tilde{A}$ uses the macroscopic discretization as $A_H$, but solves the cell problem with periodic boundary conditions in space as well as time.
	$\tilde{A}_H$ is defined via
	\begin{align*}
		\nabla \Psi_H \cdot \tilde{A}_H(t_n,x_K) \nabla \Phi_H &:= \dfrac{1}{\sigma|I_{\delta}|} \Bigl(\dfrac{\theta}{2} \int_{I_{\delta, K}} \! \nabla \Psi_H(x_K)\cdot \a0n(t_n,x) \nabla \phi^{\epsilon}(t_n,x) dx \\ &\quad  +  \theta \sum_{ k= 1}^{n-1}\int_{I_{\delta,K}} \! \nabla \Psi_H(x_K) \cdot \a0n(t_n + s_k,x) \nabla \phi^{\epsilon}(t_n + s_k,x)dx   \\ & \qquad  + \dfrac{\theta}{2} \int_{I_{\delta, K}} \! \nabla  \Psi_H(x_K)\cdot \a0n(t_n + s_{N_{cell}},x) \nabla\phi^{\epsilon}(t_n + s_{N_{cell}},x)dx \Bigr),
	\end{align*}
	where $\phi^{\epsilon}$ solves \eqref{cell_anfangs}.
	Note that $\tilde{A}_H$ includes the approximation of the temporal integral, but in contrast to $A_{H,h}$ solves the microscopic cell problems exactly.
	In the following, we write $\Phi,\Psi $ instead of $\Phi_H,\Psi_H$ for simplicity and omit the variables in the integrals for readability.
	The central result of this section is
	\begin{Satz}
		\label{eHMM_Fehlerschranke}
		If $a \in C([0,T]\times\Omega, C^2((0,1),C^1(Y) )$, then it holds for any $\alpha>0$
		\begin{align*}
			e(\text{HMM})\leq C \Bigl( \Bigl(\frac{\epsilon}{\delta}\Bigr)^{1/2} 	 	+ \dfrac{\epsilon}{\sigma^{1/2}} + \dfrac{\theta}{\epsilon^2}\dfrac{\sqrt{\sigma}}{\epsilon} +   \frac{h^{3-\alpha}}{\epsilon^3} + \frac{h}{\epsilon} \,  \Bigr),
		\end{align*}
		where $\sigma $ and $\delta$ are the time and cell size, respectively, of the cell problem \eqref{cell_anfangs}.
	\end{Satz}
	The first two terms arise from the oversampling as well as the change of the temporal and spatial boundary conditions  and are already present in \cite{ming2007analysis}. The other terms come from the discretization of the cell problems, where the leading terms orders are $h/\epsilon$ and $\theta/\epsilon^2$. We refer to  Remark~\ref{rem:h-symmetric} for a discussion of the spatial order. Linear convergence in time is expected due to the choice of implicit Euler for time stepping. To sum up, when choosing $\delta$ and $\sigma$ in practice, one has to balance the oversampling and the microdiscretization error in $e(HMM)$.
	Further, note that no additional terms $\delta$ and $\sigma$ appear as it is the case in \cite{ming2007analysis}. The reason is that we fix the macroscopic scales in the coefficient (so-called macroscopic collocation of the HMM, cf. \cite{abdulle2012heterogeneous}). 
	
	For the proof we use the triangle inequality and our auxiliary matrices $\tilde{A}$ and $\tilde{A}_H$ via
	\begin{align}
		\label{EqeMM}
		|\!| A_0 - \Ahh|\!|	\leq |\!| A_0 - \tilde{A}|\!| + |\!|   \tilde{A}  - A_H|\!| + |\!| A_H - \tilde{A}_H|\!| + |\!|  \tilde{A}_H - \Ahh|\!|.
	\end{align}
	In the following, these terms will be estimated in four steps. 
	
	\paragraph*{1.Step: Estimate $|\!| A_0 - \tilde{A}|\!|$. } 
	This error is caused by the wrong cell time and cell size. \\
	To even consider the difference of $A_0$ and $\tilde{A}$, we still need an alternative representation of $ \nabla \Psi \cdot A_0(t_n,x_K) \nabla \Phi$. 	Let $l = \lfloor \sigma/  \epsilon^2 \rfloor$, $\kappa = \lfloor \delta/\epsilon\rfloor$ and  $	\tilde{\mathcal{Q}}_{n,K}:= I_{\kappa\epsilon,K} \times (t_n,t_n + l	\epsilon^2)$.
	Since $\chin$ is the solution of problem (\ref{cell}), it follows that
	\begin{align*}
		\avint_{\tilde{\mathcal{Q}}_{n,K}} \nabla \Psi \cdot \a0n 			\nabla \hat{\Phi}^{\epsilon} &= \avint_{\tilde{\mathcal{Q}}_{n,K}} \nabla \Psi			\cdot \a0n \nabla 	\Phi + \nabla \Psi \cdot 				\a0n D_y \chin\cdot \nabla\Phi \\ &= \nabla  		\Psi  \cdot \Bigl(\avint_{\tilde{Q}_n}\a0n( \operatorname{Id}_d +  D_y 			\chin)\Bigr)  \nabla \Phi \\ &=  \nabla  \Psi			\cdot A_0(t_n, x_K)  \nabla \Phi.
	\end{align*}
	With that, we can handle the first step of estimating $e($HMM$)$. 
	\begin{Proposition}
		\label{Propositionstep1}
		There exists a constant $C$ such that
		\begin{align*}
			|\!| (\tilde{A} - A_0)|\!| \leq C\Bigl(\dfrac{\epsilon}{\delta} + 						\dfrac{\epsilon^2}{\sigma} \Bigr),
		\end{align*}
	\end{Proposition}
	\begin{proof}
		From the calculation above it follows that
		\begin{align*}
			\! \!|\nabla \Psi \cdot (A_0-\tilde{A})(t_n,x_K) \nabla \Phi| &= \Bigl|\dfrac{1}{|l\epsilon^2||I_{			\kappa\epsilon}|}\int_{I_{\kappa\epsilon,K}\times (t_n,t_n + l\epsilon^2) } \nabla \Psi \cdot \a0n \nabla \phi^{\epsilon}_{\#}  \\ & \qquad -\dfrac{1}{|\sigma||I_{				\delta}|} \int_{I_{\delta,K}\times (t_n,t_n + \sigma) } \nabla \Psi \cdot \a0n \nabla \phi^{\epsilon}_{\#}  \Bigr| \\ &= 
			\Bigl|\Bigl(\dfrac{1}{|\sigma||I_{\delta}|} - \dfrac{1}{|l\epsilon^2||I_{\kappa			\epsilon}|}\Bigr) \int_{I_{\kappa\epsilon,K}\times (t_n,t_n + l\epsilon^2)}\nabla \Psi \cdot \a0n \nabla \phi^{\epsilon}_{\#} \\ & \qquad+  			\dfrac{1}{|\sigma||I_\delta|}\int_{I_{\delta,K} \setminus I_{\kappa\epsilon,K}		\times (t_n,t_n + l\epsilon^2)}\nabla \Psi \cdot \a0n \nabla \phi^{\epsilon}_{\#} \\ & \qquad  +\dfrac{1}{|\sigma||I_{\delta}|}	\int_{I_{\delta,K} \times (t_n,t_n + (\sigma - l	\epsilon^2))}\nabla \Psi \cdot \a0n \nabla \phi^{\epsilon}_{\#} \Bigr|
			\\ &= \Bigl| 
			\underbrace{\Bigl(1-  \dfrac{|I_{\kappa\epsilon}| \cdot l \epsilon^2}{\sigma |I_{\delta}|}  		\Bigr)}_{=:G_1} \avint_{I_{\kappa\epsilon,K}\times (t_n,t_n + l\epsilon^2)}\nabla \Psi \cdot \a0n \nabla \phi^{\epsilon}_{\#} \\* & \qquad+  	\underbrace{\Bigl( \dfrac{(|I_\delta|- |I_{\kappa\epsilon}|)(l\epsilon^2)}{|\sigma| |I_			\delta|} \Bigr)}_{=:G_2} 	\avint_{I_{\delta,K} \setminus I_{\kappa\epsilon,K}\times (t_n,t_n + 		l		\epsilon^2)}\nabla \Psi \cdot \a0n \nabla \phi^{\epsilon}_{\#}\\* & \qquad  +\underbrace{\Bigl( \dfrac{|I_\delta|(\sigma -l\epsilon^2)}{|\sigma| |I_\delta|} \Bigr)}_{=:G_3}			\avint_{I_{\delta,K} \times (t_n,t_n + (\sigma - l					\epsilon^2))}\nabla \Psi\cdot \a0n \nabla \phi^{\epsilon}_{\#} \Bigr|.
		\end{align*}	
		We now consider the terms one by one. We use that  $ \delta - \epsilon \leq \kappa \epsilon \leq \delta$ and  $ \sigma - \epsilon^2 \leq \kappa \epsilon^2 \leq \sigma$ to estimate
		\begin{align*}
			G_1 &= \Bigl( 1 - \dfrac{|I_{\kappa\epsilon}| \cdot l \epsilon^2}{\sigma |I_{\delta}|}  		\Bigr)  \leq \Bigl( 1 - 				\dfrac{(\delta^d - \epsilon^d)(\sigma - \epsilon^2)}{\sigma |I_\delta|}\Bigr)			 \\ &\leq \Bigl( 						\frac{\epsilon}{\delta} + \frac{\epsilon^2}{\sigma} - \dfrac{\epsilon^{d+2}}{\sigma			\delta^d} \Bigr) \leq 					\Big(\frac{\epsilon}{\delta} + \frac{\epsilon^2}{\sigma} \Big) .
		\end{align*}
		Using the same arguments, the estimate follows for $G_2$ as well
		\begin{align*}
			G_2 &= \Bigl( \dfrac{(|I_\delta|- |I_{\kappa\epsilon}|)(l\epsilon^2)}{|\sigma| |I_			\delta|} \Bigr)	\leq 
			\Bigl(\dfrac{l\epsilon^2}{\sigma}  - \dfrac{(\delta^d - \epsilon^d)(\sigma - 				\epsilon^2)}{\sigma |I_\delta|}\Bigr)  \\ &\leq \Bigl( \frac{\epsilon}{\delta} + 						\frac{\epsilon^2}{\sigma} - \dfrac{\epsilon^{d+2}}{\sigma\delta^d} \Bigr)					 \leq \Big(\frac{\epsilon}{\delta} 		+ \frac{\epsilon^2}{\sigma} \Big) .
		\end{align*}
		To estimate $G_3$ we use that $ l\dfrac{\epsilon^2}{\sigma} \geq 1 - \dfrac{\epsilon^2}{\sigma}$ and obtain
		\begin{align*}
			G_3 = \Bigl( \dfrac{|I_\delta|(\sigma -l\epsilon^2)}{|\sigma| |I_\delta|} \Bigr)			\leq \Bigl( 1 - l \frac{\epsilon^2}{\sigma } \Bigr) \leq 					\dfrac{\epsilon^2}{\sigma}.
		\end{align*}
		Altogether, after substitution and from the boundedness of $a^{\epsilon}$  it follows that  
		\begin{equation*}
			|\nabla \Psi \cdot (\tilde{A} - A_0)(t_n,x_K) \nabla \Phi|\leq C\Bigl( \dfrac{\epsilon}{\delta} 		+ \dfrac{\epsilon^2}{\sigma} \Bigr) |\nabla \Phi| |\nabla \Psi|.\qedhere
		\end{equation*}
	\end{proof}	
	
	\paragraph*{2.Step: Estimate  $|\!|  \tilde{A} - A_H |\!|$ }
	This error can be described as the error of using wrong boundary values and wrong time conditions. \\
	Define $\zeta^{\epsilon} : = \phi^{\epsilon} -\phi^{\epsilon}_{\#}$. Then $\zeta^{\epsilon}$ satisfies the following differential equation
	\begin{align}
		\label{theta}
		\begin{cases}
			\dfrac{\partial \zeta^{\epsilon} }{\partial t} - \,\nabla \cdot(\a0n 			\nabla\zeta^{\epsilon})&= 0  \qquad \qquad \quad \; \text{ in }  \mathcal{Q}_{n,K}\\
			\zeta^{\epsilon}  &= - \epsilon \chin \nabla \Phi \quad \text{ on } \partial 				I_{\delta,K} \times 	(t_n,t_n + \sigma) 		\\
			\zeta^{\epsilon} |_{t = t_n} &= - \epsilon \chin \nabla \Phi .
		\end{cases}
	\end{align}
	We derive an estimate for $\zeta^{\epsilon}$ in the following. This in turn provides a bound for the term $\phi^{\epsilon} - \phi^{\epsilon}_{\#}$, which appears in the final calculation of the error of $ A_H $ with respect to $ \tilde{A}$. 
	\begin{Lemma}
		\label{Abschätzung_theta}
		There exists a constant $C$ independent of $\epsilon, \delta, \sigma$ such that
		\begin{align*}
			|\!| \nabla\zeta
			^{\epsilon}|\!|_{L^{2}(\mathcal{Q}_{n,K})} \leq C \biggl(\Bigl(\dfrac{\epsilon}{\delta}\Bigr)^{1/2} + \dfrac{\epsilon}{\sigma^{1/2}}\biggr) |\!| \nabla \Phi|\!|_{L^2(\mathcal{Q}_{n,K})}
		\end{align*}
		for all $\Phi \in V_H$. 
	\end{Lemma}
	The proof can be found in \cite[Lemma 3.3]{ming2007analysis}.
	We can now complete the second step as well.
	\begin{Proposition}
		\label{Propositionstep2}
		It holds that
		\begin{align*}
			|\!| (\tilde{A} - A_H) |\!| \leq C\Bigl(\Bigl( \dfrac{\epsilon}{\delta}\Bigr)^{1/2} +\dfrac{\epsilon}{\sigma^{1/2}}\Bigr).
		\end{align*}
	\end{Proposition}
	\begin{proof}
		Using Lemma \ref{Abschätzung_theta} we obtain
		\begin{align*}
			\avint_{\mathcal{Q}_{n,K}} \nabla \Psi \cdot \a0n	\nabla 	(\phi^{\epsilon}_{\#} - \phi^{\epsilon})  &\leq \dfrac{\Lambda}{|\mathcal{Q}_{n,K}|} |\!| \nabla \Psi		|\!|_{L^2(\mathcal{Q}_{n,K})} |\!| \nabla \phi^{\epsilon}_{\#} - \nabla \phi^{\epsilon}|\!|_{L^2(\mathcal{Q}_{n,K})} 
			\\ &\leq C\Bigl(  \Bigl(\dfrac{\epsilon}{\delta}\Bigr)^{1/2} +\dfrac{\epsilon}				{\sigma^{1/2}} \Bigr) \dfrac{1}{|\mathcal{Q}_{n,K}|}|\!|\nabla \Phi|\!|_{L^2(\mathcal{Q}_{n,K})} |\!| \nabla \Psi |\!|_{L^2(\mathcal{Q}_{n,K})} \\ &\leq C\Bigl( 				\Bigl(\dfrac{\epsilon}{\delta}\Bigr)^{1/2} +\dfrac{\epsilon}{\sigma^{1/2}}\Bigr)|				\nabla \Psi (x_K)||\nabla \Phi(x_K)|,
		\end{align*}
		where we used again that $\nabla \Phi$ is constant. 
	\end{proof}

	\paragraph*{3.Step: Estimate  $|\!| A_H - \tilde{A}_H|\!|$ }	
	This error can essentially be described as a quadrature error.
	\begin{Proposition}
		\label{Propositionstep3}
		For $ \Phi, \Psi \in V_H$ we obtain
		\begin{align*}
			|\nabla \Psi \cdot (A_H - \tilde{A}_H) \nabla \Phi| \leq C  \theta^2 |\nabla \Psi| |\nabla \Phi|.
		\end{align*}
	\end{Proposition}
	The proof follows directly from the quadrature order of the trapezoidal rule. In general, if use a quadrature rule of order $q$ in the definition of $A_{H,h}$ (and consequently, for $\tilde{A}_H$), this error will be bounded by $ \theta^q$.
	
	\paragraph*{4.Step: Estimate  $|\!|\tilde{A}_H - \Ahh |\!|$ }	
	This term describes the error from the microscopic discretization.
	A direct calculation shows that
	\begin{align*}
		\phi^{\epsilon}(t,x) = \Phi_H + \eta (t,x) \cdot \nabla \Phi_H(x_K)
	\end{align*}
	where $ \eta = (\eta^ 1, \dots , \eta^d ) \in \bigl(L^2((t_n,t_n + \sigma),H^1_{0}(I_{\delta,K} )) \cap H^1_{0}((t_n,t_n + \sigma),H^{-1}_{0}(I_{\delta,K} ))\bigr)^d$ satisfies
	\begin{align*}
		\int_{I_{\delta,K}} \partial_t \eta^i (t,x) z(x) + \nabla z(x) \cdot \a0n \bigl( e^i + \nabla \eta^i(t,x)\bigr)dx = 0
	\end{align*}
	 for all $z\in H^1_0(I_{\delta, K})$ and $\eta(0,x) = 0$.
	Analogously,
	\begin{align*}
		\phi_{h,k}^{\epsilon}(x) = \Phi_H(x) + \eta_{h,k}(x_K) \cdot \nabla \Phi_H(x)
	\end{align*}
	where  $ \eta_{h,k} = (\eta_{h,k}^ 1, \dots , \eta_{h,k}^d )\in (\Vmicro(I_{\delta,K}) )^d$ satisfies
	\begin{align*}
		\int_{I_{\delta,K}} \overline{\partial}_t \eta^i_{h,k}(x) z_h(x) + \nabla z_h(x) \cdot 	\a0n(t_k,x) \bigl(e^i +\nabla\eta^i_{h,k}(x)\bigr)dx = 0
	\end{align*}
	for all $z_h\in \Vmicro(I_{\delta,K})$ and $\eta_{h,0}^i=0$ for all $i=1,\ldots d$.
	$\eta$ and $\eta_{h,k}$ obviously depend on $\epsilon, \delta$ and $\sigma$.	
	
To better investigate this dependence, we re-scale $\eta$ in the following lemma.
	\begin{Lemma}
		\label{TrafoLemma}
		Define $\xi^i \in L^2((0, \frac{\sigma}{\epsilon^2}), H_0^1(I_{\delta/\epsilon}))$ via $\xi^i(s,y) = \frac{1}{\epsilon} \eta^i(t_n + \epsilon^2s, x_K + \epsilon y)$. Then $\xi^i$ solves
		\begin{align*}
			\begin{cases}
				\int_{I_{\delta/\epsilon}} \partial_s \xi^i \tilde{z} + \nabla \tilde{z} \cdot 	\tilde{a}^{\epsilon}_{n,K} \bigl(e^i +\nabla\xi^i\bigr)dx &= 0	\qquad \forall \tilde{z} \in H_0^1(I_{\delta/\epsilon})	\\
				\xi^{\epsilon} |_{t = 0} &= 0 ,
			\end{cases}
		\end{align*}
		where $\tilde{a}^{\epsilon}_{n,K}(s,y) = a(t_n,x_K,\frac{t_n}{\epsilon^2}+ s,\frac{x_K}{\epsilon} + y)$.
	\end{Lemma}
	\begin{proof}
		Define 
		\begin{align*}
			x_K^{\epsilon,\delta}: I_{\delta/\epsilon} &\to I_{\delta,K}: y \to x_K + \epsilon y \\			
			t_n^{\epsilon,\sigma}:  (0,\frac{\sigma}{\epsilon^2})&\to(t_n, t_n + \sigma) : t \to t_n + \epsilon^2 s 
		\end{align*}
		Using the transformation and chain rule, we obtain

		\begin{align*}
			&\!\!\!\int_{I_{\delta,K}} \partial_t \eta^i (t,x) z(x) + \nabla z(x) \cdot \a0n \bigl( e^i + \nabla \eta^i(t,x)\bigr)dx  \\
			 &= \int_{I_{\delta/\epsilon}} \frac{1}{\epsilon}\partial_s  \eta^i(t_n^{\epsilon,\delta}(s), x_K^{\epsilon,\delta}(y))  \frac{1}{\epsilon}\tilde{z}(x_K^{\epsilon,\delta}(y))\\
			 &\qquad + \frac{1}{\epsilon} \nabla  \tilde{z}(x_K^{\epsilon,\delta}(y)) \cdot 	\tilde{a}^{\epsilon}_{n,K}(s,y) \bigl(e^i + \frac{1}{\epsilon}\nabla \eta^i(t_n^{\epsilon,\delta}(s), x_K^{\epsilon,\delta}(y))\bigr)dy. \qedhere
		\end{align*}
	\end{proof}
Based on the rescaling of $\eta$, we estimate the error between $\eta$ and $\eta_{h,k}$, which is the key ingredient for the bound of $|\!|\tilde{A}_H - \Ahh |\!|$.
	
	\begin{Proposition}
		\label{Propositionstep4}
		Assume that  $a \in C^{2}((t_n, t_n + \sigma),C^2(I_{K,\delta}))$. Then, it holds for any $\alpha >0$
		\begin{align*}
			|\!| A_{H,h} - \tilde{A}_H|\!| \leq C \Bigl(\dfrac{\theta}{\epsilon^2}\dfrac{\sqrt{\sigma}}{\epsilon} +   \frac{h^{3-\alpha}}{\epsilon^3} +  \frac{h}{\epsilon} \Bigr) .
		\end{align*}
	\end{Proposition}
	\begin{proof}
		Let $\Phi_H, \Psi_H \in V_H$. We obtain with the boundedness of $a_{n,K}^\epsilon$ and the fact that $\nabla \Psi_H$ is piece-wise constant 
		\begin{align*}
			&\hspace*{-2ex}| \nabla \Psi_H \cdot A_{H,h} \nabla \Phi_H - \nabla \Psi_H \cdot \tilde{A}_H \nabla \Phi_H|\\ 
			&= \frac{1}{|\mathcal{Q}_{n,K}|} \Bigl( \frac{\theta}{2} \int_{I_{K,\delta}} \nabla \Psi_H\cdot \a0n(t_n,x) \nabla \bigl(\phi_{h,0}^{\epsilon}(x) - \phi^{\epsilon}(t_n,x)\bigr)  \\
			& \qquad\qquad + \theta \sum_{k = 1}^{N_{cell}-1 } \int_{I_{K,\delta}} \nabla \Psi_H(x)\cdot \a0n(t_n + s_k)  \nabla \bigl(\phi_{h,k}^{\epsilon}(x) - \phi^{\epsilon}(t_n + s_k ,x)\bigr) \\
			& \qquad\qquad +  \frac{\theta}{2} \int_{I_{K,\delta}} \nabla \Psi_H \cdot \a0n(t_n + s_{N_{cell}},x) \nabla \bigl(\phi_{h,n}^{\epsilon}(x) - \phi^{\epsilon}(t_n,x)\bigr) \Bigr) \\
			&\leq C \frac{1}{|\mathcal{Q}_{n,K}|} \Bigl( \sum^{N_{cell}}_{k = 1} \theta\|\nabla \Psi_H\|_0^2\Bigr)^{1/2}\Bigl( \sum^{N_{cell}}_{k = 1} \theta|\!| \nabla (\phi_{h,k}^{\epsilon} - \phi^{\epsilon})(t_n, \cdot) |\!|_0^2\Bigr)^{1/2}  \\ &\leq 	C \frac{1}{|\mathcal{Q}_{n,K}|} \Bigl( \sum^{N_{cell}}_{k = 1} \theta\|\nabla \Psi_H\|_0^2\Bigr)^{1/2}  \Bigl( \sum^{N_{cell}}_{k = 1} \theta|\!| \nabla (\eta^i_{h,k} - \eta^i(t_k, \cdot) )|\!|_0^2\|\nabla \Phi_H\|_0^2\Bigr)^{1/2} 	\\	
			&\leq 	C \frac{1}{|\mathcal{Q}_{n,K}|} \Bigl( \sum^{N_{cell}}_{k = 1} \theta\|\nabla \Psi_H\|_0^2\Bigr)^{1/2}  \Bigl( \sum^{N_{cell}}_{k = 1} \theta \Big[   \theta^2  \|\partial_{tt} \eta^i\|^2_{L^{2}((t_n, t_n + \sigma),L^2(I_{K,\delta}))}  \\
			& \qquad\qquad\qquad\qquad \qquad \qquad + h^{6-2\alpha} \|\partial_t \eta^i(t_k)\|^2_{H^2(I_{K,\delta})}  + h^2 \| \eta^i(t_k)\|^2_{H^2(I_{K,\delta})}\Big]\|\nabla \Phi_H\|_0^2\Bigr)^{1/2}\\	
			&\leq 	C \frac{1}{|\mathcal{Q}_{n,K}|} \Bigl( \sum^{N_{cell}}_{k = 1} \theta\|\nabla \Psi_H\|_0^2\Bigr)^{1/2}  \Bigl( \Bigl[ \sigma   \theta^2  \|\partial_{tt} \eta^i\|^2_{L^{2}((t_n, t_n + \sigma),L^2(I_{K,\delta}))} \\
			& \qquad\qquad\qquad\qquad\qquad\qquad\qquad  + h^{6-2\alpha} \|\partial_t \eta^i\|^2_{L^{2}((t_n, t_n + \sigma),H^2(I_{K,\delta}))} \\*
			& \qquad \qquad\qquad\qquad \qquad \qquad \qquad + h^2 \| \eta^i\|^2_{L^{2}((t_n, t_n + \sigma),H^2(I_{K,\delta}))}\Big]\|\nabla \Phi_H\|_0^2\Bigr)^{1/2}
		\end{align*}	
		where we used Theorem \ref{Convergenz_Theorem_parabolsisch} in the last step. We now employ Lemma \ref{TrafoLemma} and regularity results for parabolic problems \cite{sammon1982} to estimate the above terms as follows
		\begin{align*}
			\| \partial_{tt}\eta^i\|_{L^{2}((t_n, t_n + \sigma),L^2(I_{K,\delta}))} & = \epsilon \frac{1}{\epsilon^4} {\epsilon^{d/2}} \epsilon
			\|\partial_{ss}\xi^i\|_{L^{2}((0,\frac{\sigma}{\epsilon^2}),L^2(I_{\delta/\epsilon}))}  \leq C_2\frac{1	}{\epsilon^3} {\epsilon^{d/2}}\sqrt{   |I_{\delta/\epsilon}|} \, \epsilon \, \sqrt{\dfrac{\sigma}{\epsilon^2}}  \\ 
			&\leq C_2\frac{1	}{\epsilon^3} \sqrt{  |I_{\delta,K}|} \sqrt{{\sigma}},\\
			\|\partial_t \eta^i\|_{{L^{2}((t_n, t_n + \sigma),H^2(I_{K,\delta}))}}  & = \epsilon \frac{1}{\epsilon^2} \frac{1}{\epsilon^2} {\epsilon^{d/2}} \epsilon
			\|\partial_{s}\xi^i\|_{{L^{2}((0, \frac{\sigma}{\epsilon^2}),H^2(I_{\delta/\epsilon}))}}   \leq C_3\frac{1	}{\epsilon^3} {\epsilon^{d/2}}\sqrt{   |I_{\delta/\epsilon}|} \,\epsilon \, \sqrt{\dfrac{\sigma}{\epsilon^2}} \\
			&\leq C_3\frac{1}{\epsilon^3} \sqrt{  |I_{\delta,K}|}\sqrt{{\sigma}}, \\ 
			\| \eta^i\|_{{L^{2}((t_n, t_n + \sigma),H^2(I_{K, \delta}))}}  & = \epsilon \frac{1}{\epsilon^2}  {\epsilon^{d/2}}\epsilon 
			 \|\xi^i\|_{{L^{2}((0, \frac{\sigma}{\epsilon^2}),H^2(I_{\delta/\epsilon}))}}  \leq C_4\frac{1	}{\epsilon} {\epsilon^{d/2}}\sqrt{|I_{\delta/\epsilon}|}\, \epsilon \,\sqrt{\dfrac{\sigma}{\epsilon^2}} \\
			&\leq C_4\frac{1	}{\epsilon} \sqrt{  |I_{\delta,K}|}\sqrt{{\sigma}}.
		\end{align*}
		Inserting these inequalities, we get
		\begin{align*}
			| \nabla \Psi_H \cdot A_{H,h} \nabla \Phi_H - \nabla \Psi_H \cdot \tilde{A}_H \nabla \Phi_H| &\leq  C \Bigl(\dfrac{\theta}{\epsilon^2}\dfrac{\sqrt{\sigma}}{\epsilon} +   \frac{h^{3-\alpha}}{\epsilon^3} +  \frac{h}{\epsilon} \Bigr) | \nabla \Phi_H(x_K)|\, |\nabla \Psi_H(x_K)|.\qedhere
		\end{align*}
		
	\end{proof}
	\begin{Bemerkung}\label{rem:h-symmetric}
		Since we are in the non-symmetric case we only get a theoretical convergence order of $h$. For the case that $A_0$ is symmetric, convergence order of $h^2$ is expected. In a different setting with non-symmetric homogenized coefficient, \cite{freese2022heterogeneous} even observed $h^2$ convergence numerically. 
		The term $h^{3-\alpha}$ occurs due to the $H^{-1}$-norm estimate and is discussed in the appendix.
	\end{Bemerkung}
	Summing up, we have proved the estimate for $e($HMM$)$.
	\begin{proof}[Proof of Theorem \ref{eHMM_Fehlerschranke}.]
		Using  \eqref{EqeMM} and Propositions \ref{Propositionstep1}, \ref{Propositionstep2}, \ref{Propositionstep3} and \ref{Propositionstep4} we obtain the  desired error bound.
	\end{proof}

	\section{Numerical Experiments \label{result}}
	In the following, numerical results calculated with the Heterogeneous Multiscale Method are presented. The  convergence rate with respect to the time step and mesh width is investigated for the macrodiscretization as well as for the microdiscretization.
	The implementation was done in Python, building on the \textit{Fenics} software library \cite{langtangen2017solving}, where Version 2019.2.0 was used for this paper.

	\subsection{Setting}
	We choose $\Omega = (0,1)$ and $T = 1$. Let the initial condition be $ u^{\epsilon}(0,x) = 0$ for all $x \in \Omega$.
	As the exact solution of the homogenized equation we choose 
	\begin{align*}
		U_{0}(t,x) = t^2(x - x^2)
	\end{align*}
	and accordingly the right side  
	\begin{align*}
		f = 2t(x-x^2) + 2A_0 t^2
	\end{align*}
	with the homogenized coefficient $A_0$.
	As in Section \ref{Errmacro}, the error between the numerical solution $U_H$ and the exact solution $U_0$ is investigated. For this purpose, we consider $|\!| U_H^N - U_0(t_N,x)|\!|_{L^2(\Omega)}$ the error at time $t_N$.

	As seen in Theorem \ref{eHMM_Fehlerschranke}, the error bound of e(HMM) depends on the terms $\sigma$, $\frac{\epsilon}{\delta}$, and $\frac{\epsilon^2}{\sigma}$, which we get from the boundary and initial values. 
	When choosing the parameters $\delta$ and $\sigma$, it is important that $H \gg \delta$ and $H \gg \sigma$. Here, as suggested by  \cite{ming2007analysis}, we choose $\delta = \epsilon^{1/3} $ and $\sigma = \epsilon^{2/3}$.

	\subsection{First example }  
	In the first example, we select the coefficient as
	\begin{align*}
		a(t,x,s,y) = 3 + \cos(2\pi y) + \cos^2(2\pi s).
	\end{align*}
	According to \cite{tan2019high},
	$A_0 \approx 3.352429824667637$. For simplicity and more efficient calculation, in $f$ the coefficient $A_0$ is replaced by this value.
	 
\begin{figure}
		\centering
		\includegraphics[width=0.45\textwidth]{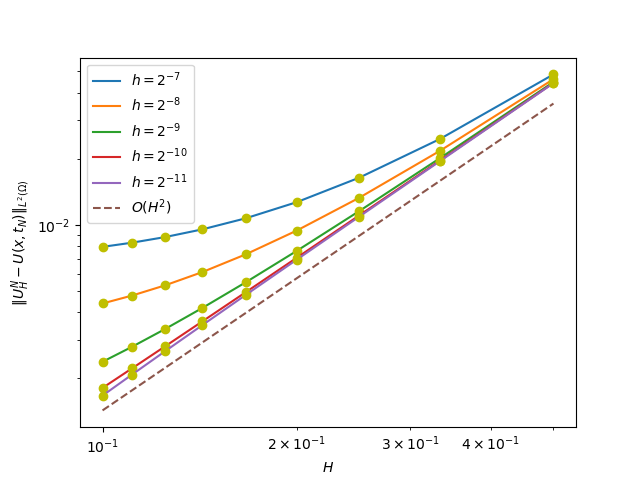}%
		\hspace{2ex}%
		\includegraphics[width=0.45\textwidth]{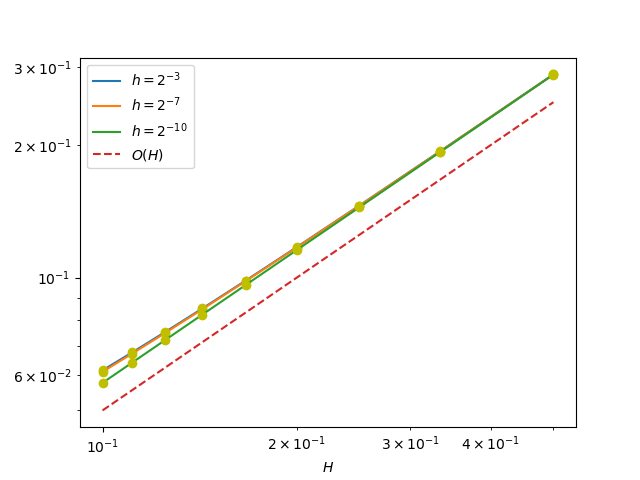}
		\caption{$L^2$-error (left) and $H^1$-error (right) with respect to  the grid sizes $H$ and $h$ of the macro and micro discretization, respectively,  at the time $t = 1$ and $\epsilon = 10^{-3}$.  }
		\label{Bsp0}
	\end{figure}

	Figure \ref{Bsp0} (left) shows the error in the $L^2$ norm over the grid width $H$ for different cell grid widths $h$. Time step sizes were fixed as $\tau = \frac{1}{15}$ and $\theta = \frac{\sigma}{15}$. It can be clearly seen that for coarser cell grid widths $h = 2^{-7} $ and $h = 2^{-8}$ the error of microdiscretization dominates and therefore we do not get convergence order $H^2$. However, for fine cell grid widths, the expected order shows up. 
	Similarly, we obtain linear convergence w.r.t. $H$ in the $H^1(\Omega)$ norm (Figure \ref{Bsp0} right). Note that for the $H^1$ norm, convergence in the macro mesh size can already be observed for relatively large cell grid widths.
	These results agree nicely with our findings in Theorems \ref{Fehler_Satz} and \ref{eHMM_Fehlerschranke}.
	\subsubsection{Second example }
\begin{figure}
		\centering
		\includegraphics[width=0.45\textwidth]{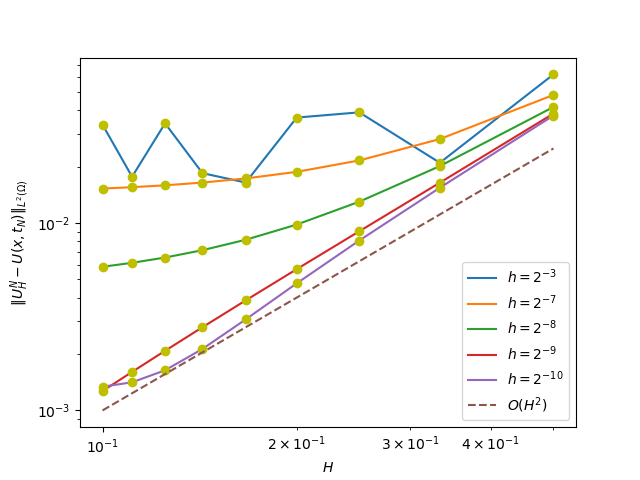}%
		\hspace{2ex}%
		\includegraphics[width=0.45\textwidth]{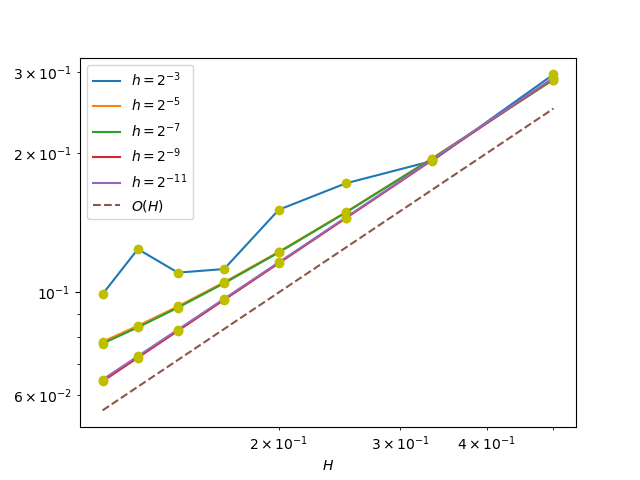}
		\caption{$L^2$-error (left) and $H^1$-error (right) with respect to  the grid sizes $H$ and $h$ of the macro and micro discretization, respectively,  at the time $t = 1$ and $\epsilon = 10^{-3}$. }
		\label{Bsp1}
	\end{figure}

	In the next example we set
	\begin{align*}
		a_1(y) = \dfrac{1}{2 - \cos(2\pi y)}
	\end{align*}
	independent of $s$. 
	A straight forward calculation shows
	\begin{align*}
		A_0 = C_0 = \dfrac{1}{2}.
	\end{align*}
	We again consider the error between the numerical solution $U_H$ and the exact solution $U_0$. 
	As in the previous example, we choose as time step sizes $\tau = \frac{1}{15}$ and $\theta = \frac{\sigma}{15}$.
	We again first consider the $L^2$-error and study its convergence w.r.t. to $H$, see Figure \ref{Bsp1} left.
	For $h = 2^{-3}$ no convergence can be seen. For fine grid widths $h = 2^{-7} $ and $h = 2^{-8}$, quadratic convergence is initially seen, but still the microerror dominates. Only for very fine cell grid widths quadratic convergence does appear. The fact that for $h = 2^{-10}$ the error flattens out is probably due to rounding errors. Figure \ref{Bsp1} (right) shows the convergence in the $H^1$ norm. Here, for sufficiently small cell grid width, the expected linear convergence is also observed. The results are again in alignment with the theory and the observations for the first example.
	
\begin{figure}
		\centering
		\includegraphics[width=0.6\textwidth]{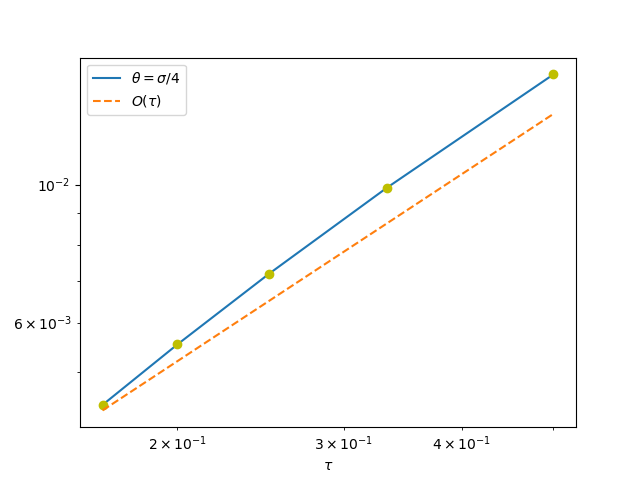}
		\caption{$L^2$-error with respect to time steps $\tau$   at $t = 1$ and $\epsilon = 10^{-3}$. }
		\label{Bsp1_L2_time}
	\end{figure}

	In addition, we investigate the error of the time discretization. For this, we fix $H = \frac{1}{10}$ and $ h= \frac{1}{1000}$. Figure \ref{Bsp1_L2_time} shows that for fixed cell time step size $ \theta = \frac{\sigma}{4}$ linear convergence in the macro time step $\tau$ can be observed. Again, this underlines our theoretically predicted results.

	\bibliographystyle{abbrv}
	\bibliography{source}
	
	\appendix
	\section{Finite element error estimates for parabolic problems}
	In this appendix, we present and prove an a priori error estimate for the finite element discretization of an (abstract) parabolic problem with time-dependent coefficient. The statement and its proof are very similar to the well known results in the literature cf., e.g., \cite{thomee2006}. We try to extract as  high spatial convergence order as possible for each term. This is crucial to obtain (almost) balanced orders between $h$ and $\epsilon$ in Proposition \ref{Propositionstep4}.
	In the following, we use $\|\cdot \|_m$ to denoted the $H^m(\Omega)$-norm for any $m\in \R$.
			\begin{Satz}
				\label{Convergenz_Theorem_parabolsisch}
				Let $\Omega \subset \mathbb{R}^d$ be a  hypercube, $ T >0$, $b:[0,T] \times H^1(\Omega) \times H^1(\Omega)$ be bounded, coercive and Lipschitz-continuous (w.r.t. the time $t$) and $f:[0,T] \to L^2(\Omega)$ be continuously differentiable. Denote by  $u$ the solution of  the following weak problem
				\begin{align}
					\label{parablic_weak}
					\begin{cases}
						\Big\langle\dfrac{\partial u}{\partial t}, v\Big\rangle_0 -  b(t,u,v ) &= (f,v)_0  \quad \forall  v \in H_0^1(\Omega) \\
						\langle u(0+), v\rangle_0 &= 0 \qquad \forall v\in  L^2(\Omega).
					\end{cases}
				\end{align} 
				Let further $u_h^n \in \Vmicro$ be the fully discrete approximation of $u$ at time $t_n = n\tau$ using implicit Euler method with time step size $\tau$ and  finite elements of order $p\geq 2$. Then we have for any $\alpha>0$
				\begin{align}
					\notag
					\Bigl(\tau \sum_{j= 1}^n \| u_h^j - u(t_j)\|_{H_0^1(\Omega)}^2 \Bigr)^{1/2} \leq
					C \Bigl( \sum^{n}_{j = 1} \tau \Big[ \tau^2  \|\partial_{tt} u\|^2_{L^{2}((0, T),L^2(\Omega))}
					 \label{convergenz_para}
					 &+ h^{6-2\alpha} \|\partial_t u(t_k)\|^2_{H^2(\Omega)} \\  & + h^2 \| u(t_k)\|^2_{H^2(\Omega)}\Big]\Bigr)^{1/2}  .
				\end{align} 
			\end{Satz}
			\begin{proof}
			We denote by $R_h$ the Ritz projection onto $\Vmicro$ with respect to $b(t,
		\cdot,\cdot)$. Note that $R_h$ depends on the time $t$, but we will omit this dependence for better readability.
				Inserting the Ritz projection of the exact solution into the discrete equation and using a standard stability estimate, cf., e.g., \cite{thomee2006}, provides
				\begin{align*}
					\Bigl(\tau \sum_{j= 1}^n \| u_h^j - R_h u(t_j)\|_{H_0^1(\Omega)}^2 \Bigr)^{1/2} \leq \Bigl(\tau \sum_{j= 1}^n \| d_n\|_{H_0^{-1}(\Omega)}^2 \Bigr)^{1/2} ,
				\end{align*}
			where 
			\begin{align*}
				d_n= \frac{1}{\tau}\int_{t_{n-1}}^{t_n}\partial_t((I-R_h)u)\, dt+\dfrac{u(t_n) - u(t_{n-1})}{\tau} - \partial_t u(t_n) .
			\end{align*}
			The first term can be bounded with the following arguments of \cite{thomee2006}: Set $e = u - R_h u$. Let $w \in H^1(\Omega)$ be arbitrary and find $z \in H_0^1(\Omega) $ such that $b(\cdot,v,z) = (w,v)_{0} $. From elliptic regularity estimates it follows $z \in H^{3-\alpha}(\Omega)$ for any $\alpha>0$ since $\Omega$ is a hypercube \cite{bourland1992}. If we choose $ v= \partial_t e$, we obtain
			\begin{align*}
				(\partial_t e, w)_{0} = b(\cdot,\partial_te, z ) = b(\cdot,\partial_te, z - v_h) + b^{\prime}(\cdot,e,z - v_h) - b^{\prime}(\cdot,e,z),
			\end{align*}
			where the second equality follows by differentiating the equation $b(\cdot,e,v_h) = 0$ for all $v_h \in \Vmicro$ with $b^{\prime}(\cdot,\cdot,\cdot)$ the
			bilinear form obtained from $b(\cdot,\cdot, \cdot)$ by differentiating the coefficients with
			respect to $t$.  We get 
			\begin{align*}
				(\partial_t e, w)_0 \leq (\|\partial_te\|_1 + \|e\|_1 ) \inf_{v_h \in \Vmicro} \|z - v_h\|_1 + \|e\|_{{-1}}\|z\|_{3-\alpha}
			\end{align*}
			By elliptic regularity and convergence results for (at least) quadratic finite elements \cite{bourland1992,brenner-scott}
			\begin{align*}
				(\partial_t e, w) \leq h^{2-\alpha}(\|\partial_te\|_1 + \|e\|_1 ) \|w\|_{1} + \|e\|_{{-1}}\|w\|_{1}.
			\end{align*}
			From \cite{bourland1992,brenner-scott} it follows for small $h$
			\begin{align*}
				\|\partial_t e\|_{-1} \leq h^{3-\alpha}(\|\partial_t u\|_{2} +  \|u\|_{2}) + h^{3-\alpha} \|u\|_{2}
				\leq Ch^{3-\alpha}(\|\partial_t u\|_{2} +\|u\|_{2}).
			\end{align*}

			For the second term in $d_n$ we use again projection estimates and Taylor expansion
			\begin{align*}
				\|\Bigl(\dfrac{u(t_n) - u(t_{n-1})}{\tau} - \partial_t u(t_n) \Bigr)\|_{H^{-1}(\Omega)}&\leq  \|\Bigl(\dfrac{u(t_n) - u(t_{n-1})}{\tau} -  \partial_t u(t_n) \Bigr)\|_{L^2(\Omega)} \\  &\leq \frac{1}{2} \tau \Bigl(  \int_{t_{n-1}}^{t_n}\|\partial_{tt}u(s)\|^2_{L^2(\Omega)} \,ds\Bigr)^{1/2}  .
			\end{align*}
			Using the triangle inequality and estimates for the Ritz projection error of the exact solution, cf. \cite{brenner-scott,thomee2006}, we obtain \eqref{convergenz_para}.
			\end{proof}
		
		We emphasize that quadratic finite elements are necessary to get the required bounds for the Ritz projection error in the $H^{-1}$-norm, cf. \cite{bourland1992,brenner-scott}. Since our domain of interest is a hypercube, we only get $H^{3-\alpha}$-regularity for any $\alpha>0$ of a solution to an elliptic problem with right-hand side in $H^1$.
	
\end{document}